\documentclass{amsart}
\usepackage{latexsym,amsmath,amssymb}

\theoremstyle{definition}
\theoremstyle{remark}

\numberwithin{equation}{section}

\begin{document}

\title[Multiplication conditional
type operators ]
{Fredholm Multiplication conditional
type operators on $L^{p}$-space }

\author{\sc Y. Estaremi}
\email{estaremi@gmail.com} 

\address{Department of mathematics, Payame Noor university , P. O. Box: 19395-3697, Tehran, Iran}


\thanks{}

\thanks{}

\subjclass[2010]{47B38, 47B37}

\keywords{Conditional expectation, multiplication operators, Fredholm,
closed range.}

\date{}

\dedicatory{}

\commby{}

\begin{abstract}
In this paper, first we investigate closed range  multiplication conditional
type operators between two $L^p$-spaces. Then we characterize Fredholm ones when the underlying measure space is non-atomic. Finally we give some examples.

\noindent {}
\end{abstract}

\maketitle

\section{ \sc\bf Introduction and Preliminaries}

Let $(X,\Sigma,\mu)$ be a complete $\sigma$-finite measure space.
For any sub-$\sigma$-finite algebra $\mathcal{A}\subseteq
 \Sigma$ with $1\leq p\leq \infty$, the $L^p$-space
$L^p(X,\mathcal{A},\mu|_{\mathcal{A}})$ is abbreviated  by
$L^p(\mathcal{A})$, and its norm is denoted by $\|.\|_p$. All
comparisons between two functions or two sets are to be
interpreted as holding up to a $\mu$-null set. The support of a
measurable function $f$ is defined as $S(f)=\{x\in X; f(x)\neq
0\}$. We denote the vector space of all equivalence classes of
almost everywhere finite valued measurable functions on $X$ by
$L^0(\Sigma)$.

\vspace*{0.3cm} For a sub-$\sigma$-finite algebra
$\mathcal{A}\subseteq\Sigma$, the conditional expectation
operator associated with $\mathcal{A}$ is the mapping
$f\rightarrow E^{\mathcal{A}}f$, defined for all non-negative function $f$
as well as for all $f\in L^p(\Sigma)$, $1\leq p\leq \infty$,
where $E^{\mathcal{A}}f$, by the Radon-Nikodym theorem, is the
unique $\mathcal{A}$-measurable function satisfying
$$\int_{A}fd\mu=\int_{A}E^{\mathcal{A}}fd\mu, \ \ \ \forall A\in \mathcal{A} .$$
As an operator on $L^{p}({\Sigma})$, $E^{\mathcal{A}}$ is an
idempotent and $E^{\mathcal{A}}(L^p(\Sigma))=L^p(\mathcal{A})$. If
there is no possibility of confusion we write $E(f)$ in place of
$E^{\mathcal{A}}(f)$. This operator will play a major role in our
work and we list here some of its useful properties:

\vspace*{0.2cm} \noindent $\bullet$ \  If $g$ is
$\mathcal{A}$-measurable, then $E(fg)=E(f)g$.

\noindent $\bullet$ \ $|E(f)|^p\leq E(|f|^p)$.

\noindent $\bullet$ \ If $f\geq 0$, then $E(f)\geq 0$; if $f>0$,
then $E(f)>0$.

\noindent $\bullet$ \ $|E(fg)|\leq
E(|f|^p)|^{\frac{1}{p}}E(|g|^{p'})|^{\frac{1}{p'}}$, where
$\frac{1}{p}+\frac{1}{p'}=1$ (H\"{o}lder inequality).

\noindent $\bullet$ \ For each $f\geq 0$, $\sigma(f)\subseteq
\sigma(E(f))$.

\vspace*{0.2cm}\noindent A detailed discussion and verification of
most of these properties may be found in \cite{rao}. We recall that an $\mathcal{A}$-atom of the measure
$\mu$ is an element $A\in\mathcal{A}$ with $\mu(A)>0$ such that
for each $F\in\mathcal{A}$, if $F\subseteq A$ then either
$\mu(F)=0$ or $\mu(F)=\mu(A)$. A measure space $(X,\Sigma,\mu)$
with no atoms is called non-atomic measure space. It is well-known
fact that every $\sigma$-finite measure space $(X,
\mathcal{A},\mu_{\mid_{\mathcal{A}}})$ can be partitioned uniquely
as $X=\left (\bigcup_{n\in\mathbb{N}}A_n\right )\cup B$, where
$\{A_n\}_{n\in\mathbb{N}}$ is a countable collection of pairwise
disjoint $\mathcal{A}$-atoms and $B$, being disjoint from each
$A_n$, is non-atomic (see \cite{z}).

 \vspace*{0.2cm} Compositions
of conditional expectation operators and multiplication operators
appear often in the study of other operators such as
multiplication operators and weighted composition operators.
Specifically, in \cite{mo}, S.-T. C. Moy characterized all
operators on $L^p$ of the form $f\rightarrow E(fg)$ for $g$ in
$L^q$ with $E(|g|)$ bounded. Eleven years later, R. G. Douglas,
\cite{dou}, analyzed positive projections on $L^{1}$ and many of
his characterizations are in terms of combinations of
multiplications and conditional expectations. More recently, P.G.
Dodds, C.B. Huijsmans and B. De Pagter, \cite{dhd}, extended these
characterizations to the setting of function ideals and vector
lattices. J. Herron presented some assertions about the operator
$EM_u$ on $L^p$ spaces in \cite{ her}. Also, some results about
multiplication conditional type operators can be found in
\cite{g, lam}. In \cite{e,ej3,ej} we investigated some classic
properties of multiplication conditional type operators
$M_wEM_u$ on $L^p$ spaces. Let $f\in L^0(\Sigma)$, then $f$ is
said to be conditionable with respect to $E$ if
$f\in\mathcal{D}(E):=\{g\in L^0(\Sigma): E(|g|)\in
L^0(\mathcal{A})\}$. Throughout this paper we take $u$ and $w$ in
$\mathcal{D}(E)$. In this paper, some necessary and sufficient
conditions for closeness of range of multiplication conditional
type operators between two $L^p$-spaces are given. Also, Fredholm ones are characterized when the underlying measure space is non-atomic.  Some results of this paper are a generalization of some results of \cite{tak}.\\

Now we give a definition of multiplication conditional type operators on $L^p$-spaces.

\vspace*{0.4cm} {\bf Definition 1.1.}
Let $(X,\Sigma,\mu)$ be a $\sigma$-finite measure space and let $\mathcal{A}$ be a
$\sigma$-subalgebra of $\Sigma$ such that $(X,\mathcal{A},\mu_{\mathcal{A}})$ is also $\sigma$-finite. Let $E$ be the corresponding conditional
expectation operator relative to $\mathcal{A}$. If $1\leq p,q<\infty$ and $u,w \in L^0(\Sigma)$ (the spaces of $\Sigma$-measurable functions on $X$) such that $uf$ is conditionable and $wE(uf)\in L^{q}(\Sigma)$ for all $f\in \mathcal{D}\subseteq L^{p}(\Sigma)$, where $\mathcal{D}$ is a linear subspace, then the corresponding multiplication conditional type operator is the linear transformation $M_wEM_u:\mathcal{D}\rightarrow L^{q}(\Sigma)$ defined by $f\rightarrow wE(uf)$.\\

 Here we recall some results of \cite{dhd} that state our results is valid for a large class of linear operators. Let $(X, \Sigma, \mu)$ be a finite measure space, then $L^{\infty}(\Sigma)\subseteq L^p(\Sigma)\subseteq L^1(\Sigma)$ and $L^p(\Sigma)$ is an order ideal of measurable functions on $(X,\Sigma,\mu)$. Thus by propositions $(3.1, 3.3, 3.6)$ of \cite{dhd} we have theorems A, B, C:\\

\vspace*{0.4cm} {\bf Theorem A.} If $T$ is a linear operator on $L^p(\Sigma)$ for which\\

(i) $Tf\in L^{\infty}(\Sigma)$ whenever $f\in L^{\infty}(\Sigma)$.\\

(ii) $\|Tf_n\|_1\rightarrow 0$ for all sequences $\{f_n\}_{n=1}^{\infty}\subseteq L^{p}(\Sigma)$ such that $|f_n|\leq g$ $(n=1,2,3,....)$ for some $g\in L^{p}(\Sigma)$ and $f_n\rightarrow 0$ a.e.,\\

(iii) $T(f.Tg)=Tf.Tg$ for all $f\in L^{\infty}(\Sigma)$ and all $g\in L^{p}(\Sigma)$,\\

then there exists a $\sigma$-subalgebra $\mathcal{A}$ of $\Sigma$ and there exists $w\in L^{p'}(\Sigma)$ such that $Tf=E^{\mathcal{A}}(wf)$ for all $f\in L^{p}(\Sigma)$.$(p^{-1}+p'^{-1}=1)$\\

 \vspace*{0.4cm} {\bf Theorem B.} For a linear operator $T:L^p(\Sigma)\rightarrow L^p(\Sigma)$ the following statement are equivalent.\\

(i) $T$ is positive and order continuous, $T^2=T$, $T1=1$ and the range $\mathcal{R}(T)$ of $T$ is a sublattice.\\

(ii) There exist a $\sigma$-subalgebra $\mathcal{A}$ of $\Sigma$ and a function $0\leq w\in L^{p'}(\Sigma)$ with $E^{\mathcal{A}}(w)=1$  such that $Tf=E^{\mathcal{A}}(wf)$ for all $f\in L^{p}(\Sigma)$.\\

\vspace*{0.4cm} {\bf Theorem C.} For a linear operator $T:L^p(\Sigma)\rightarrow L^p(\Sigma)$ the following statement are equivalent.\\

(i) $T$ is a positive and order continuous projection onto a sublattice such that $T1$ is strictly positive.\\

(ii) There exist a $\sigma$-subalgebra $\mathcal{A}$ of $\Sigma$, $0\leq w\in L^{p'}(\Sigma)$  and a strictly positive function $k\in L^1(\Sigma)$ with $E^{\mathcal{A}}(wk)=1$  such that $Tf=kE^{\mathcal{A}}(wf)$ for all $f\in L^{p}(\Sigma)$. Moreover, if we choose $k$ such that $E^{\mathcal{A}}(k)=1$, then both $w$ and $k$ are uniquely determined by $T$. \\

\section{ \sc\bf CLOSED RANGE AND FREDHOLM WEIGHTED CONDITIONAL TYPE OPERATORS }
In this section first we describe closed range multiplication conditional type operators $M_wEM_u$ between two $L^p$-spaces. Let $1\leq p,q<\infty$ and $f\in L^p$, then it is easily seen that 
$$\|M_wEM_u(f)\|_q=\|EM_v(f)\|_q,$$
where $v=u(E(|w|^q))^{\frac{1}{q}}$. Thus without loss of generality we can consider the operator $EM_v$
instead of $M_wEM_u$, in our discussion about closedness of range. Also, we recall that for an operator $T$ on a Banach space $X$, $\mathcal{N}(T)=\{x\in X:T(x)=0\}$ and $\mathcal{R}(T)=\{T(x): x\in X\}$ are called null space and range of $T$, respectively.

\vspace*{0.3cm} {\bf Theorem 2.1.} Let $1<p<\infty$ and let $p'$
be conjugate component to $p$. Then\\

(a) If the operator $T=EM_u$ from $L^{p}(\Sigma)$ into itself is injective and has closed range, then there
exists $\delta>0$ such that
$v=(E(|u|^{p'}))^{\frac{1}{p'}}\geq\delta$ a.e, On $S$, where
$S=\{x\in X:v(x)\neq0\}$.\\

(b) If $\sigma(E(u))=\sigma(E(|u|^{p'}))$ and there exists
$\delta>0$ such that $E(u)\geq\delta$ a.e. On
$S=\sigma(E(|u|^{p'}))=\{x\in X:E(|u|^{p'})(x)\neq0\}$, then the
operator $T=EM_u$ has closed range on $L^{p}(\Sigma)$.

\vspace*{0.3cm} {\bf Proof.} (a) Let $f\in L^{p}(\Sigma)$. Then
\begin{align*} 
\|Tf\|^{p}_{p}&=\int_{X}|E(uf)|^{p}d\mu\\
&\leq\int_{X}(E(|u|^{p'}))^{\frac{p}{p'}}E(|f|^{p})d\mu\\
&=\int_{X}v^p|f|^{p}d\mu=\|M_{v}f\|^{p}_{p}.
\end{align*}

Since $T$ is injective and closed range, then there exists $\delta>0$ such that for $f\in L^{p}(\Sigma)$, $\|Tf\|_{p}\geq\delta\|f\|_{p}$. Thus
\begin{align*}
\|M_{v}f\|_{L^p(S)}&=\|M_{v}f\|_{L^p(X)}\\
&\geq\|Tf\|_{p}\\
&\geq\delta\|f\|_{L^p(X)}\\
&\geq\delta\|f\|_{L^p(S)}
\end{align*}
and so $\|M_{v}f\|_{L^p(S)}\geq\delta\|f\|_{L^p(S)}$, for all $f\in L^{p}(\Sigma)$.
This mean's that $M_{v}$ has closed range on $L^{p}(X)$. Thus there exists $\beta>0$ such that $v\geq\beta$ a.e. On $S$.

(b) Let $f_n,g\in L^{p}(\Sigma)$ such that
$\|E(uf_n)-g\|_p\rightarrow0$, when $n\rightarrow\infty$. Since
$E(u)\geq\delta$ a.e. On $S$, then $\frac{1}{E(u)}\leq
\frac{1}{\delta}$ a.e. On $S$. This implies that
$\frac{g}{E(u)}\chi_{S}\in L^p(S)$ and
$E(u\frac{g}{E(u)}\chi_{S})=g\in L^p(X, \mathcal{A}, \mu)$. Hence
\begin{align*}
\|E(uf_n)-E(u\frac{g}{E(u)}\chi_{S})\|^p_p&=\int_{X}|E(uf_n)-E(u\frac{g}{E(u)})\chi_{S}|^pd\mu\\
&=\int_{S}|E(uf_n)-E(u\frac{g}{E(u)})\chi_{S}|^pd\mu\\
&=\int_{S}|E(uf_n)-g|^pd\mu\\
&\leq\|E(uf_n)-g\|^p_p\rightarrow0,
\end{align*}

when $n\rightarrow\infty$. So the operator $EM_u$ has closed range
on $L^{p}(\Sigma)$.

\vspace*{0.3cm} {\bf Theorem 2.2.} Let $1<q<p<\infty$ and let $p',
q'$ be conjugate component to $p$ and $q$ respectively. Then

(a) If the operator $T=EM_u$ from $L^{p}(\Sigma)$ into
$L^{q}(\Sigma)$ is injective and has closed range, then
\\
\begin{enumerate}
\item  $v=0$ a.e. On $B$ and the set $\{n\in \mathbb{N}:v(A_{n})\neq0\}$ is finite.\\
\item  $M_{v}$ from $L^{p}(\Sigma)$ into $L^{q}(\Sigma)$ has finite rank.\\
\end{enumerate}
Where $v=(E|u|^{q'})^{\frac{1}{q'}}$ and $S=\{x\in X:v(x)\neq0\}$.

(b) Let\\
\begin{enumerate}
\item The operator $T=EM_u$ has closed range.\\

\item The operator $T=EM_u$ has finite rank.\\

\item $v=0$ a.e. On $B$ and the set $N_{v}=\{n\in \mathbb{N}:v(A_{n})\neq0\}$ is finite.\\

\item $E(u)=0$ a.e. On $B$ and the set $N_{N(u)}=\{n\in \mathbb{N}:E(u)(A_{n})\neq0\}$ is finite.\\
\end{enumerate}
Then $$(3)\rightarrow (2)\rightarrow (1)\rightarrow (4).$$

 \vspace*{0.3cm} {\bf Proof.} (a)
 Let $f\in L^{p}(\Sigma)$. Then
 \begin{align*}
\|Tf\|^{q}_{q}&=\int_{X}|E(uf)|^{q}d\mu\\
&\leq\int_{X}(E(|u|^{q'}))^{\frac{q}{q'}}E(|f|^{q})d\mu\\
&=\int_{X}v^q|f|^{q}d\mu\\
&=\|M_{v}f\|^{q}_{q}.
\end{align*}

Since $T$ is injective and closed range, then there exists $\delta>0$ such that for $f\in L^{p}(\Sigma)$, $\|Tf\|_{q}\geq\delta\|f\|_{p}$. Thus
$\|M_{v}f\|_{q}\geq\|Tf\|_{p}\geq\delta\|f\|_{p}$
and so $\|M_{v}f\|_{q}\geq\delta\|f\|_{p}$, for all $f\in L^{p}(\Sigma)$.
This mean's that $M_{v}$ from $L^{p}(\Sigma)$ into
$L^{q}(\Sigma)$ has closed range. Thus by \cite{tak} we have $v=0$ a.e. On $B$ and the
 set $\{n\in \mathbb{N}:v(A_{n})\neq0\}$ is finite. Also, $M_{v}$ from $L^{p}(\Sigma)$ into $L^{q}(\Sigma)$ has finite
 rank.\\

(b) If $\mu(S)=0$ then $EM_u$ is the zero operator. So we assume
that $\mu(S)>0$. $(3)\rightarrow (2)$. If $(3)$ holds, then
$S=\cup_{n\in N_{v}}A_n=\cup^{k}_{i=1}A_{n_i}$ for some integer
$k>0$. Hence $$EM_u(L^p(X, \Sigma, \mu))\subseteq L^p(S,
\mathcal{A}, \mu),$$ Since for any $f\in L^p(X, \Sigma, \mu)$,
$\sigma(E(uf))\subseteq S$. This implies that $EM_u$ has finite
rank.\\

$(2)\rightarrow (1)$ is trivial.\\

$(1)\rightarrow (4)$. Suppose that $EM_u$ has closed range.First
we show that $E(u)=0$ a.e. On $B$. Assume on the contrary that
$\mu(\{x\in B: E(u)(x)\neq0\})>0$. Then we have $\mu(\{x\in B:
E(u)(x)>\delta\})>0$ for some $\delta>0$. Set $G=\{x\in
B:E(u)(x)>\delta\}$ and define a function $v$ on $G$ by
$v(x)=\frac{1}{E(u)(x)}$ for $x\in G$. For every $f\in L^p(G)$ and
$g\in L^q(G)$,
$$M_vEM_u(f)=v(E(u)|_G)f=f, \ \ \ \ M_{E(u)|_G}M_v(g)=g.$$
Thus $M_v$ is the inverse operator of $EM_u|_{L^p(G)}=M_{E(u)|_G}$
and $EM_u|_{L^p(G)}=M_{E(u)|_G}$ is a bounded operator from
$L^p(G)$ into $L^q(G)$ that has closed rang.\\

For any $E\in \mathcal{A}_G=\{A\cap G: A\in \mathcal{A}\}$ with
$\mu(E)<\infty$. put $f=\frac{1}{E(u)(x)}\chi_{E}(x)$, then $f\in
L^p(G)$. Moreover, $M_{E(u)}f=\chi_{E}$ and so $\chi_E\in
M_{E(u)}(L^p(G))$. Hence, the the range  $M_{E(u)}(L^p(G))$
contains all linear combinations of such $\chi_E$'s. Thus
$M_{E(u)}(L^p(G))=L^q(G)$, since $M_{E(u)}$ has closed range and
all linear combinations of such $\chi_E$'s are dense in $L^q(G)$.
This implies that $M_v$ maps $L^q(G)$ into $L^p(G)$, that is $M_v$
is bounded from $L^q(G)$ into $L^p(G)$. Since $G$ is non-atomic by
theorem 1.4 of \cite{tak} we have $v=0$ a.e on $G$. But this is a
contradiction.\\

 Now, we show that $N_{E(u)}$ is finite. Since
$\mu(S)>0$, it follows that $N_{E(u)}\neq \emptyset$. Put
$w(x)=\frac{1}{E(u)(x)}$ for $x\in S$ and consider the operator
$M_w$. put $f=\frac{1}{E(u)(x)}\chi_{A_n}(x)$, then $f\in L^p(G)$,
by the same method of preceding paragraph we see that $M_w$ maps
$L^q(S)$ into $L^p(S)$ that is $M_w$ is bounded from $L^q(S)$ into
$L^p(S)$. So by theorem 1.4 of \cite{tak} we have
$$b=\sup_{n\in N_{E(u)}}\frac{1}{|E(u)(A_n)|^r\mu(A_n)}=\sup_{n\in
N_{E(u)}}\frac{|w(A_n)|^r}{\mu(A_n)}<\infty,$$ where
$\frac{1}{p}+\frac{1}{r}=\frac{1}{q}$. Since $N_{E(u)}\neq
\emptyset$, then $b>0$ and $|E(u)(A_n)|^r\mu(A_n)\geq\frac{1}{b}$.
While Theorem 1.3 of \cite{tak} says $E(u)\in L^r(X, \mathcal{A},
\mu)$. So,

$$\Sigma_{n\in N_{E(u)}}\frac{1}{b}\leq\|E(u)\|^r_r<\infty.$$
This implies that $N_{E(u)}$ is finite.\\

\vspace*{0.3cm} {\bf Theorem 2.3.} Let $1< p< q<\infty$ and let
$p', q'$ be conjugate component to $p$ and $q$ respectively.\\

(a) If
the operator $T=EM_u$ from $L^{p}(\Sigma)$ into $L^{q}(\Sigma)$ is
injective and has closed range, then
\\
\begin{enumerate}
\item The set $\{n\in \mathbb{N}:v(A_{n})\neq0\}$ is finite.\\
\item $M_{v}$ from $L^{p}(\Sigma)$ into $L^{q}(\Sigma)$ has finite rank.\\
\end{enumerate}
Where $v=(E|u|^{q'})^{\frac{1}{q'}}$ and
$S=\{x\in X:v(x)\neq0\}$.\\

(b)  Let\\
\begin{enumerate}
\item The operator $T=EM_u$ has closed range.\\

\item The operator $T=EM_u$ has finite rank.\\

\item The set $N_{v}=\{n\in \mathbb{N}:v(A_{n})\neq0\}$ is finite.\\

\item The set $N_{N(u)}=\{n\in \mathbb{N}:E(u)(A_{n})\neq0\}$ is finite.\\
\end{enumerate}
Then $$(3)\rightarrow (2)\rightarrow (1)\rightarrow (4).$$

\vspace*{0.3cm} {\bf Proof.} (a)
 Let $f\in L^{p}(\Sigma)$. Then
$\|Tf\|^{q}_{q}\leq\|M_{v}f\|^{q}_{q}.$

Since $T$ is injective and closed range, then there exists $\delta>0$ such that for $f\in L^{p}(\Sigma)$, $\|Tf\|_{q}\geq\delta\|f\|_{p}$. Thus
$\|M_{v}f\|_{q}\geq\|Tf\|_{p}\geq\delta\|f\|_{p}$
and so $\|M_{v}f\|_{q}\geq\delta\|f\|_{p}$, for all $f\in L^{p}(\Sigma)$.
This mean's that $M_{v}$ from $L^{p}(\Sigma)$ into
$L^{q}(\Sigma)$ has closed range. Thus by \cite{tak} the set $\{n\in \mathbb{N}:v(A_{n})\neq0\}$ is
 finite. Also, $M_{v}$ from $L^{p}(\Sigma)$ into $L^{q}(\Sigma)$ has finite
 rank.\\

 (b) Theorem 2.3 of \cite{ej} tells us that $v=0$ a.e on $B$ and so $E(u)=0$ a.e on $B$. By
 the same method that we used in last theorem, it is easy to see
 that $(3)\rightarrow (2)\rightarrow (1)$. Now, we show that $(1)\rightarrow
 (4)$. Suppose that $N_{N(u)}\neq \emptyset$. If we put
 $S=\sigma(E(u))$, then we can write $S=\cup_{n\in N_{E(u)}}A_n$.
 Define a function $w$ on $S$ by $w(x)=\frac{1}{E(u)(x)}$ for
 $x\in S$. By the same method that is used in the proof of last
 $L^p(S)$. Hence by Theorem 1.3 of \cite{tak} we have $w\in
 L^s(\mathcal{A})$, where $\frac{1}{q}+\frac{1}{s}=\frac{1}{p}$.
 While Theorem 1.4 of \cite{tak} says that $b=\sup_{n\in
 N_{E(u)}}\frac{|E(u)(A_n)|^s}{\mu(A_n)}<\infty$. Since $N_{N(u)}\neq
 \emptyset$ implies $b>0$ and since
 $|w(A_n)|^s\mu(A_n)=\frac{\mu(A_n)}{|E(u)(A_n)|^s}\geq\frac{1}{b}$
 for all $n\in N_{N(u)}$,
it follows that
$$\Sigma_{n\in N_{E(u)}}\frac{1}{b}\leq\|w\|^s_s<\infty$$
this implies that $N_{E(u)}$ is finite.\\

In the sequel we consider the function $u$ is $\mathcal{A}$-measurable. In this case $E(u)=u$ and $EM_u=M_u\mid_{L^p(\mathcal{A})}$. Therefore we get the following results.\\

\vspace*{0.3cm} {\bf Corollary 2.4.} Let $1<p<\infty$ and let $p'$
be conjugate component to $p$. If
$u\in L^0(\mathcal{A})$ and $T=EM_u$ is injective on $L^p(\Sigma)$. Then the operator $T=EM_u$ has closed range if and only if there exists $\delta>0$ such that
$|u|\geq\delta$ a.e. On $S$. Where
$S=\{x\in X:u(x)\neq0\}$.\\

\vspace*{0.3cm} {\bf Corollary  2.5.} Let $1<q<p<\infty$ and let $p',
q'$ be conjugate component to $p$ and $q$ respectively. If
$u\in L^0(\mathcal{A})$ and $T=EM_u$ is injective from $L^{p}(\Sigma)$ into $L^{q}(\Sigma)$. Then the followings are equivalent:

\begin{enumerate}
\item The operator $M_u$ has closed range.\\

\item The operator $M_u$ has finite rank.\\

\item $u=0$ a.e. On $B$ and the set $N_{u}=\{n\in \mathbb{N}:u(A_{n})\neq0\}$ is finite.\\
\end{enumerate}

\vspace*{0.3cm} {\bf Corollary  2.6.} Let $1< p< q<\infty$ and let
$p', q'$ be conjugate component to $p$ and $q$ respectively. If
$u\in L^0(\mathcal{A})$ and $T=EM_u$ is injective from $L^{p}(\Sigma)$ into $L^{q}(\Sigma)$, then the followings are equivalent:\\

\begin{enumerate}
\item The operator $EM_u$ has closed range.\\

\item The operator $EM_u$ has finite rank.\\

\item The set $N_{u}=\{n\in \mathbb{N}:u(A_{n})\neq0\}$ is finite.
\end{enumerate}

If the multiplication conditional type operator $EM_u$ is bounded from $L^p(\Sigma)$ onto $L^p(\mathcal{A})$, then $\mu(Z(E(|u|^{p'}))=\{x\in X:E(|u|^{p'})(x)=0\})=0$. Suppose that $F\subseteq Z(E(|u|^{p'}))$ with $F\in \mathcal{A}$ and $\mu(F)<\infty$. Then $\chi_{F}\in L^p(\mathcal{A})=\mathcal{R}(T)$ and there exists $f\in L^p(\Sigma)$ such that $Tf=\chi_{F}$. So by conditional-type H\"{o}lder inequality we have
\begin{align*}
\mu(F)&=\int_F|E(u.f)|^pd\mu\\
&\leq\int_F(E(|u|^{p'})^{\frac{p}{p'}}.|f|^pd\mu=0.
\end{align*}
Hence $\mu(Z(E(|u|^{p'})))=0$.

 In the next theorem we characterize Fredholm multiplication conditional type operators on $L^p$-spaces when the underlying measure space is non-atomic.\\

\vspace*{0.3cm} {\bf Theorem 2.7.} Let $1< p<\infty$ and $(X,\Sigma, \mu)$ be a non-atomic measure space. If
$T=EM_u$ is a bounded operator from $L^{p}(\Sigma)$ into $L^{p}(\mathcal{A})$, then $T$ is Fredholm if and only if $T$ is invertible.\\
\vspace*{0.3cm} {\bf Proof.} Let $T$ be Fredholm, then $T$ has closed range. First we show that $T$ is surjective. Suppose on the contrary. Let $f_0\in L^p(\mathcal{A})\setminus \mathcal{R}(T)$. Then there exists a bounded linear functional $L_{g_0}$ on $L^p(\mathcal{A})$ for some $g_0\in L^{p'}(\mathcal{A})$ ($p^{-1}+p'^{-1}=1$), which is defined as
$$L_{g_0}(f)=\int_Xf\bar{g_0}d\mu, \ \ \ \ \ f\in L^p(\mathcal{A}),$$
such that
$$L_{g_0}(f_0)=1,  \ \ \ \ \ L_{g_0}(\mathcal{R}(T))=0.$$
Then there exists a positive constant $\delta$ such that
$$\mu(E_{\delta}=\{x\in X: f_0(x)g_0(x)>\delta\})>0.$$
Since the underlying measure space is non-atomic, then we can find a disjoint sequences $\{E_n\}_n$ such that $E_n\subseteq E_{\delta}$ and $0<\mu(E_n)<\infty$. Let $g_n=g_0.\chi_{E_n}$. Clearly $g_n\in L^{p'}(\mathcal{A})$ and for every $f\in L^p(\mathcal{A})$ we have
\begin{align*}
\int_X\bar{f}T^{\ast}(g_n)d\mu&=\int_X\bar{f}\bar{u}E(g_n)d\mu\\
&=\int_Xg_0E(\bar{u}\bar{f}.\chi_{E_n})d\mu\\
&=L_{g_0}(Tf)=0.
\end{align*}
This implies that $g_n\in \mathcal{N}(T^{\ast})$ and so $\mathcal{N}(T^{\ast})$ is infinite dimensional. Therefore the codimension of $\mathcal{R}(T)$ is not finite. This is a contradiction, therefore $T$ is surjective. Let $0\neq f\in L^p(\Sigma)$ such that $T(f)=0$. Since $S(f)\subseteq S(E(|f|))$, $\mu(S(f))>0$ and $S(E(|f|))$ is $\mathcal{A}$-measurable, then there exists $A\in \mathcal{A}$ with $0<\mu(A)<\infty$ and $\mu(S(f)\cap A)>0$. Also since the underlying measure space is non-atomic, we can find a disjoint sequence of $\mathcal{A}$-measurable subsets $\{A_n\}_n$ of $A$ with $\mu(S(f)\cap A_n)>0$. Clearly $f_n=f.\chi_{S(f)\cap A_n}\in L^p(\Sigma)$ and $T(f_n)=\chi_{A_n}.T(f)=0$. This means $\mathcal{N}(T)$ is infinite dimensional, that is a contradiction. Thus $T$ should be injective and so is invertible. The converse is obvious.\\

In the sequel we present some examples of conditional expectations
and corresponding multiplication conditional type operators.\\

\vspace*{0.3cm} {\bf Examples.} (a) Let $(X_1,\Sigma_1, \mu_1)$ and $(X_2,\Sigma_2, \mu_2)$ be two
$\sigma$-finite measure spaces and $X=X_1\times X_2$,
$\Sigma=\Sigma_1\times \Sigma_2$ and $\mu=\mu_1\times \mu_2$. Put
$\mathcal{A}=\{A\times X_2:A\in \Sigma_1\}$. Then $\mathcal{A}$ is
a sub-$\sigma$-algebra of $\Sigma$ and for all $f\in \mathcal{D}(E^{\mathcal{A}})$ we have
$$E^{\mathcal{A}}(f)(x_1,x_2)=\int_{X_2}f(x_1,y)d\mu_2(y) \ \ \
\mu-a.e.$$ on $X$.\\

Mpreover, if $(X,\Sigma, \mu)$ is a finite measure space and
$k:X\times X\rightarrow \mathbb{C}$ is a $\Sigma\otimes
\Sigma$-measurable function such that
 the operator
$T:L^p(\Sigma)\rightarrow L^q(\Sigma)$ defined by
$$Tf(x)=\int_{X}k(x,y)f(y)d\mu, \ \ \ \ \ f\in L^p(\Sigma),$$
is bounded. Then $T$ is a
weighted conditional type operator. Since $kf$ is a $\Sigma\otimes \Sigma$-measurable function, when $f\in
L^p(\Sigma)$. Then by taking $u:=k$ and $f'(x,y)=f(y)$, we get
that
\begin{align*}
E^{\mathcal{A}}(uf)(x,y)&=E^{\mathcal{A}}(uf')(x,y)\\
&=\int_{X}u(x,y)f'(x,y)d\mu(y)\\
&=\int_{X}u(x,y)f(y)d\mu(y)\\
&=Tf(x).
\end{align*}
Hence $T$ is a weighted conditional type operator. In particular the convolution operator
$$w*u(y)=\int_G w(y-x)u(x)dm(x)$$
on an abelian group $G$ with Haar measure $m$, can be viewed as an
integral operator.\\

(b) This example deals with a
weighted conditional expectation operator which is in the form of
an integral operator. Let $X=(0,1]\times (0,1]$, $d\mu=dxdy$,
$\Sigma$  the  Lebesgue measurable subsets of $X$ and let
$\mathcal{A}$ be the $\sigma$-algebra generated by the family of
the sets of the form $A\times (0,1]$ where $A$ is a Lebesgue
measurable subset of $(0,1]$. For $f$ in
$L^2((0,1]^2)$, we have $(E^{\mathcal{A}}f)(x, y)=\int_0^1f(x,t)dt$.\\

Let $X=(0,\infty)\times (0,\infty)$ and $\mathcal{A}$ be the
$\sigma$-algebra generated by the family of the sets of the form
$A\times (0,\infty)$ where $A$ is a Lebesgue measurable subset of
$(0,\infty)$. Then
$$Tf(x,y)=(E^{\mathcal{A}}M_u(f))(x,y)=\int_0^{\infty}u(x,t)f(x,t)dt.$$
Hence for every function $f:(0,\infty)\rightarrow \mathbb{R}$ we
can define the function $f'$ on $X$ as $f'(x,y)=f(y)$. So,
\begin{align*}
Tf(x,y)&=(E^{\mathcal{A}}M_u(f))(x,y)\\
&=(E^{\mathcal{A}}M_u(f'))(x,y)\\
&=\int_0^{\infty}u(x,t)f(t)dt\\
&=T'(f)(x).\\
\end{align*}

This implies that $T$ is an integral transform and specially by
taking $u(x,y)=e^{-xy}$ we obtain one of the most important
classical integral transforms that is widely used in analysis,
namely the Laplace integral transform. We refer to \cite{lu} for
some applications of integral transforms, especially Laplace
integral transforms.\\


\begin{thebibliography}{99}
\bibitem{dhd} P.G. Dodds, C.B. Huijsmans and B. De Pagter,
characterizations of conditional expectation-type operators,
Pacific J. Math. {\bf 141}(1) (1990), 55-77.

\bibitem{e} Y. Estaremi, Essential norm of weighted conditional type operators on $L^p$-spaces, positivity. {\bf 18} (2014), 41-52.




\bibitem{ej3} Y. Estaremi and M.R. Jabbarzadeh, Compact weighted Lambert type operators on $L^p$-spaces, J. Math. Anal. Appl. {\bf 420} (2014)118-123.


\bibitem{ej} Y. Estaremi and M.R. Jabbarzadeh, Weighted lambert type operators on
$L^{p}$-spaces, Oper. Matrices {\bf 1} (2013), 101-116.



\bibitem{dou}
 R. G. Douglas, Contractive projections on an $L\sb{1}$ space,
 Pacific J. Math. {\bf 15} (1965), 443-462.

\bibitem{g}
J. J. Grobler and B. de Pagter, Operators representable as
multiplication-conditional expectation operators, J. Operator
Theory {\bf 48} (2002), 15-40.

\bibitem{her}
J. Herron, Weighted conditional expectation operators, Oper.
Matrices {\bf 1} (2011), 107-118.



\bibitem{lam}
A. Lambert, $L^p$ multipliers and nested sigma-algebras, Oper.
Theory Adv. Appl. {\bf 104} (1998),  147-153.


\bibitem{lu} Y. Luchko, Integral transforms of the Mellin convolution type and
their generating operators, Integral Trans. and special
Func. {\bf 11} (2008), 809-851.

\bibitem{mo}
Shu-Teh Chen, Moy,  Characterizations of conditional expectation
as a transformation on function spaces,
 Pacific J. Math. {\bf 4} (1954), 47-63

\bibitem{rao}
M. M. Rao, Conditional measure and applications, Marcel Dekker,
New York, 1993.

\bibitem{sha}
J. H. Shapiro, The essential norm of a composition operator,
Analls of Math. {\bf 125} (1987), 375-404.

\bibitem{tak}
H. Takagi and K. Yokouchi, Multiplication and composition
operators between two $L^p$-spaces, Contemporary Math. {\bf
232}(1999), 321-338.



\bibitem{z}
A. C. Zaanen, Integration, 2nd ed., North-Holland, Amsterdam,
1967.


\end{thebibliography}
\end{document}